\documentclass[onecolumn]{article}

\usepackage{graphicx} % handles graphics and figures
\usepackage{titling}   %Fancy title editing commands
\usepackage{fullpage}  %Change margins to be better
\usepackage{mathtools}
\usepackage{bbold}   %Get fancy double struck math notation for sets
\usepackage{courier}   %Have code written out nicely
\usepackage{placeins}   %Force figures to stay in the correct place
\usepackage{color}

\catcode`\^^M=10      %  Makes blank lines meaningless to TeX,
\definecolor{darkgreen}{rgb}{0.0, 0.4, 0.2}
\definecolor{darkgreen}{rgb}{0.0, 0.4, 0.2}

\begin{document} %

%====================================================
\setlength{\droptitle}{-60pt}  %Adjust title position -- needs titling packag
\title{Transcription Methods for Trajectory Optimization \\
\vspace{4pt}
\Large{ A beginners tutorial}}
\author{Matthew P. Kelly\\
	Cornell University\\
	\texttt{mpk72@cornell.edu}}
\date{February 18, 2015}  %Comment this line and the program will use today's date
\maketitle

%====================================================
%\renewcommand{\abstractname}{Renamed Abstract}% Rename the abstract

\begin{abstract}

This report is an introduction to transcription methods for trajectory optimization techniques. The first few sections describe the two classes of transcription methods (shooting \& simultaneous) that are used to convert the trajectory optimization problem into a general constrained optimization form. The middle of the report discusses a few extensions to the basic methods, including how to deal with hybrid systems (such as walking robots). The final section goes over a variety of implementation details.

\end{abstract}

%=============================
\section{Optimal Control Overview}

There are three types of algorithms for solving optimal control problems\cite{Diehl2009}:
\begin{itemize}
\setlength{\itemsep}{0pt}
	\item{\bf Dynamic Programming: } Solve Hamilton-Jacobi-Bellman Equations over the entire state space.
	\item{\bf Indirect Methods: } Transcribe problem then find where the slope of the objective is zero.
	\item{\bf Direct Methods: } Transcribe problem then find the minimum of the objective function.
\end{itemize}

Dynamic programming is an excellent solution to the optimal control problem for unconstrained low-dimensional systems, but it does not scale well to high-dimensional systems, since it requires a discretization of the full state space. Indirect methods tend to be numerically unstable and are difficult to implement and initialize \cite{Betts2010}. For the remainder of this paper we will restrict focus to direct methods for transcribing and solving the optimal control problem. The solution to an optimal control problem via transcription scales well to high-dimensional systems, but yields a single trajectory through state and control space, rather than a global policy like dynamic programming.

\subsection{Trajectory Optimization Problem}
\label{sec:TrajectoryOptimizationProblem}

A trajectory optimization problem seeks to find the a trajectory for some dynamical system that satisfies some set of constraints while minimizing some cost functional. Below, I've laid out a general framework for a trajectory optimization problem.\\ \\

\begin{align}

\text{\bf Optimal Trajectory: }& \quad
\{ \mathbf{x}^*(t), \mathbf{u}^*(t) \}
\label{eqn:optimalTrajectory}\\

\text{\bf System Dynamics: }& \quad
\dot{\mathbf{x}} = \mathbf{f}(t,\mathbf{x},\mathbf{u})
\label{eqn:systemDynamics}\\

\text{\bf Constraints: }& \quad
\mathbf{c}_{\text{min}} < \mathbf{c}(t,\mathbf{x},\mathbf{u}) < \mathbf{c}_{\text{max}}
\label{eqn:constraints}\\

\text{\bf Boundary Conditions: }& \quad
\mathbf{b}_{\text{min}} < \mathbf{b}(t_0,\mathbf{x}_0,t_f,\mathbf{x}_f) < \mathbf{b}_{\text{max}}
\label{eqn:boundaryConditions}\\

\text{\bf Cost Functional: }& \quad
J = \phi(t_0,\mathbf{x}_0,t_f,\mathbf{x}_f) + \int_{t_0}^{t_f} g(t,\mathbf{x},\mathbf{u}) \; dt
\label{eqn:costFunctional}

\end{align}

One interesting thing to point out is the difference between a {\em state} $\mathbf{x}$ and a {\em control} $\mathbf{u}$. A state variable is a variable that is differentiated in the dynamics equation, where as a control variable only appears algebraically in the dynamics equation \cite{Betts2010}. In some cases there might also be unknown {\em parameters} (not shown) that are time-invariant variables that appear algebraically in the dynamics equation.

\subsection{Non-linear Programming}

Transcription methods for solving an optimal control problem work by converting a continuous problem (Section \ref{sec:TrajectoryOptimizationProblem}) into a non-linear programming problem (\ref{eqn:NonLinearProgramming}). Once in this form, the problem can be passed to a commercial solver, such as SNOPT\cite{Gill2005}, IPOPT\cite{Wachter2006}, or FMINCON \cite{MatlabOptimizationToolbox2014}.

\begin{align}
	\underset{\mathbf{z} \in \mathbb{R}^n}{\text{minimize}}  \quad  & \bar{J}(\mathbf{z}) \\
	\text{subject to: }  & \mathbf{l} \leq \left( 
			\begin{array}{c}
			\mathbf{z} \\
			\mathbf{\bar{c}}(\mathbf{z}) \\
			A \mathbf{z}
			\end{array}
			\right) \leq \mathbf{u}
			\label{eqn:NonLinearProgramming}	
\end{align}

There are many transcription algorithms that make this conversion, but they can all be divided up into two broad classes: {\bf shooting methods} and {\bf simultaneous methods}. The difference is based on how each method enforces the constraint on the system's dynamics. Shooting methods use a simulation to explicitly enforce the system dynamics. Simultaneous methods enforce the dynamics at a series of points along the trajectory.

%~~~~~~~~~~~~~~~~~~~~~~~~~~~~~~~~~~~~~~~~~~~~~~~~~~~~~~~~~~~~~~~~~~~~~~~~~~~~~~~~~%
%~~~~~~~~~~~~~~~~~~~~~~~~~~~~~~~~~~~~~~~~~~~~~~~~~~~~~~~~~~~~~~~~~~~~~~~~~~~~~~~~~%
%~~~~~~~~~~~~~~~~~~~~~~~~~~~~~~~~~~~~~~~~~~~~~~~~~~~~~~~~~~~~~~~~~~~~~~~~~~~~~~~~~%

\section{Shooting Methods}

Single-shooting is probably the simplest method for transcribing an optimal control problem. Consider the problem of trying to hit a target with a cannon. You have two decision variables (firing angle and mass of powder) and one constraint (trajectory passes through target). The dynamics are simple (projectile motion) and the cost function is the mass of powder. The single shooting method is similar to what a person might accomplish with experiments. You make a guess at the angle and amount of powder, and then fire the canon. If you shoot over the target, then perhaps you would reduce the mass of powder on the next test. By repeating this method, you would eventually be able to hit the target, while using as little powder as possible. Single shooting works the same way, just replacing the experiments with simulations.

\par In a more general case of single shooting, with a continuous control input, you would choose some arbitrary function to approximate the input. A few common choices are zero-order-hold, piecewise linear, piecewise cubic, or orthogonal polynomials. If the control is modeled with a piecewise function, then you must take care to align the discontinuities of the control with the integration steps in the simulation.

\par Single shooting works well enough for simple problems, but it will almost certainly fail on problems that are more complicated. This is because the relationship between the decision variables and the objective and constraint functions is not well approximated by the linear (or quadratic) model that the non-linear programming solver uses. 

\par Multiple shooting works by breaking up a trajectory into some number of segments, and using single shooting to solve for each segment. As the segments get shorter, the relationship between the decision variables and the objective function and constraints becomes more linear. In multiple shooting, the end of one segment will not necessarily match up with the start of the next. This difference is known as a defect, and it is added to the constraint vector. Adding all of the segments will increase the number of decision variables (the start of each segment) and the number of constraints (defects). Although it might seem that this would make the low-level optimization problem harder, it actually turns out to make it easier.

\par
\begin{figure}
\centering
\includegraphics[width=\textwidth]{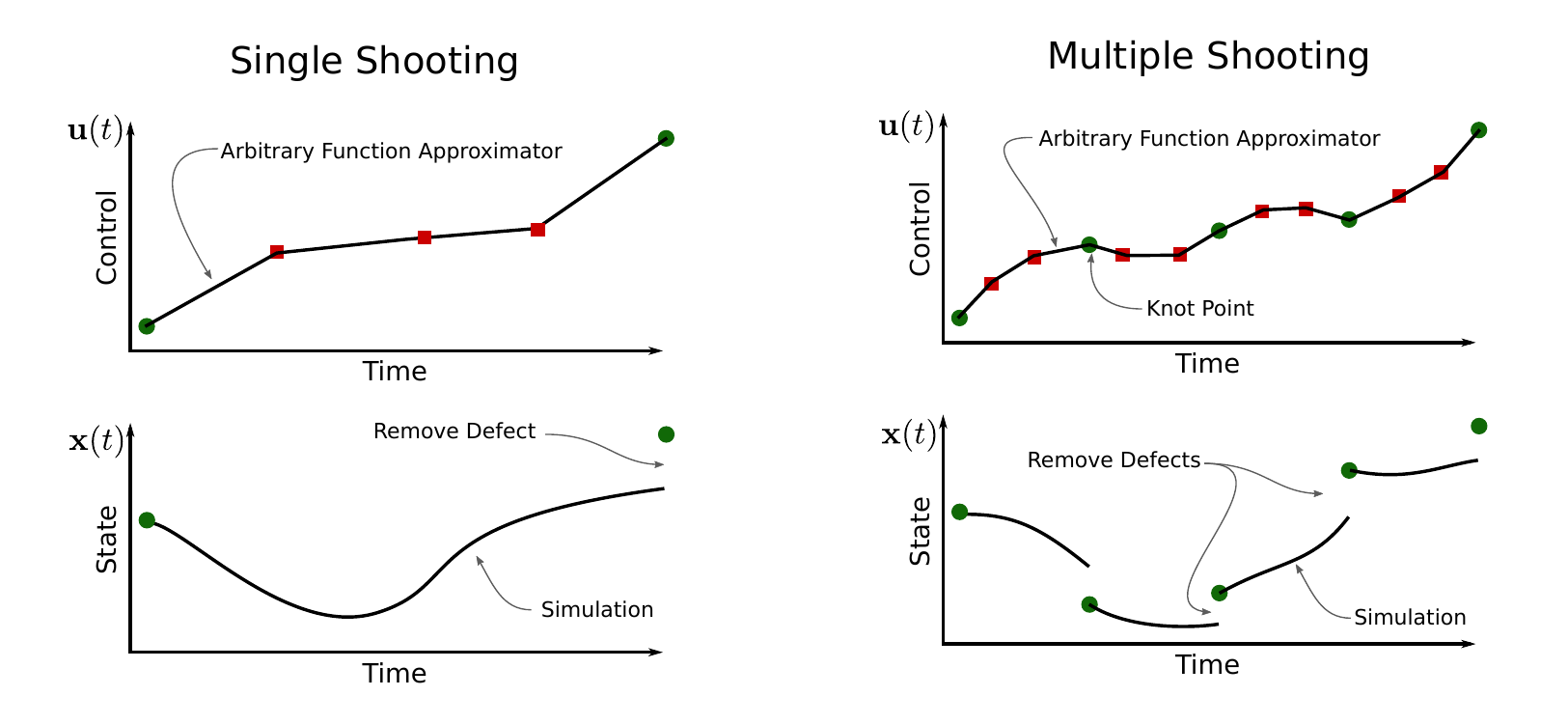}
\caption{Single Shooting vs Multiple Shooting. In both methods the state trajectory is stored as the result of a simulation. Notice that multiple shooting is just like a series of single shooting methods, with a defect constraint added to make the trajectory continuous. Multiple shooting results in a higher dimensional non-linear program, but it is sparse and more linear than the program produced by single shooting.}
\label{fig:Shooting_Methods}
\end{figure}

 Figure \ref{fig:Shooting_Methods} shows a cartoon comparing single shooting and multiple shooting.

%~~~~~~~~~~~~~~~~~~~~~~~~~~~~~~~~~~~~~~~~~~~~~~~~~~~~~~~~~~~~~~~~~~~~~~~~~~~~~~~~~%
%~~~~~~~~~~~~~~~~~~~~~~~~~~~~~~~~~~~~~~~~~~~~~~~~~~~~~~~~~~~~~~~~~~~~~~~~~~~~~~~~~%
%~~~~~~~~~~~~~~~~~~~~~~~~~~~~~~~~~~~~~~~~~~~~~~~~~~~~~~~~~~~~~~~~~~~~~~~~~~~~~~~~~%

\section{Simultaneous Methods}

There are a wide variety of simultaneous methods. The key difference between simultaneous methods and shooting methods is that simultaneous methods directly represent the state trajectory using decision variables, and then satisfy the dynamics constraint only at special points in the trajectory.

\subsection{Integral vs Differential Form}

\begin{figure}
\centering
\includegraphics[width=\textwidth]{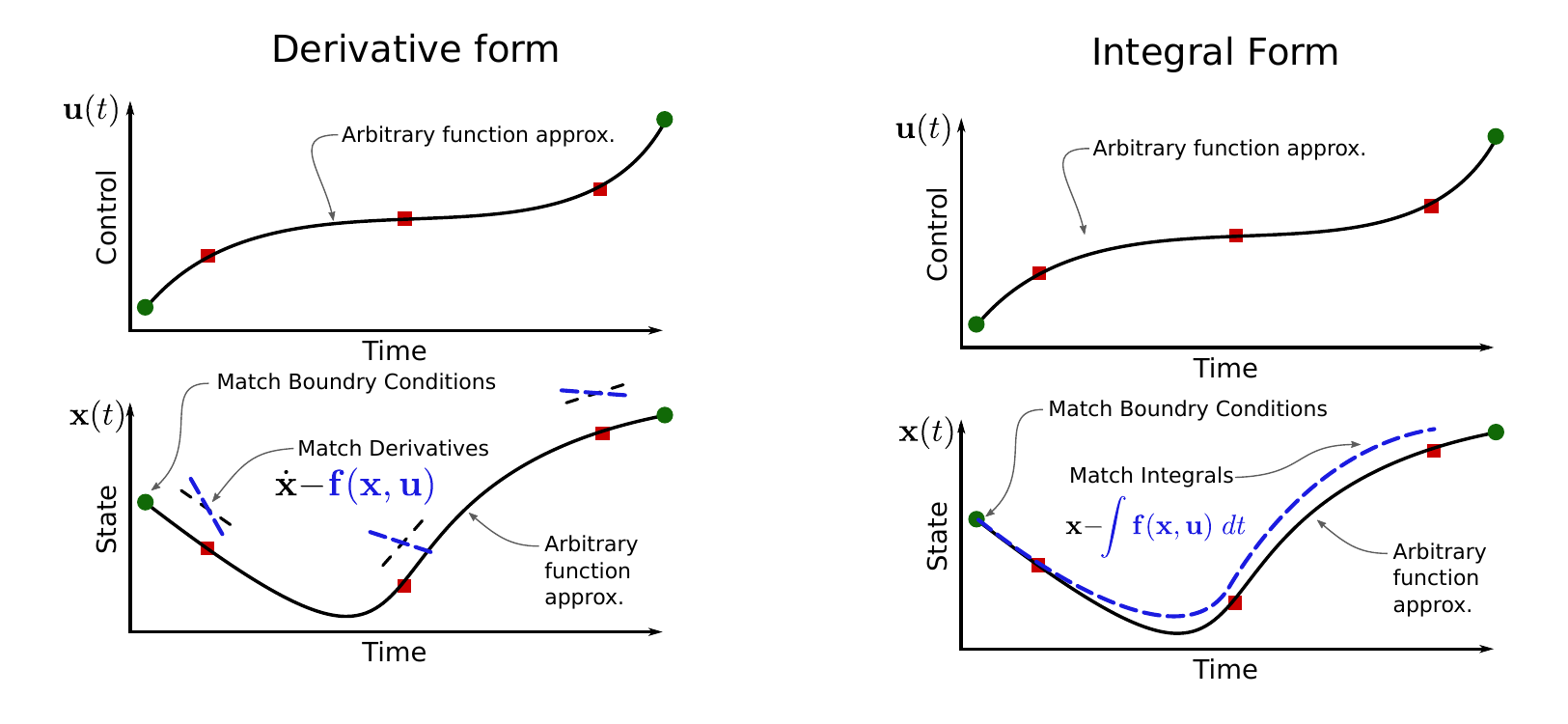}
\caption{Integral vs. derivative form of the optimal control problem.}
\label{fig:Integral_vs_Derivative}
\end{figure}

For any trajectory optimization problem, there are two different ways to write the dynamics constraint: {\em derivative} and {\em integral}. The derivative method states that the derivative of the state with respect to time must be equal to the dynamics function ($\dot{\mathbf{x}} = \mathbf{f}(\mathbf{x}, \mathbf{u})$). The integral method states that the state trajectory must match the integral of the dynamics with respect to time ($\mathbf{x} = \int \mathbf{f}(\mathbf{x}, \mathbf{u}) \; dt$). Notice that shooting methods are a type of integral method. Figure \ref{fig:Integral_vs_Derivative} shows a cartoon of the difference, and more information can be found in the paper by Fran\c{c}olin {\em et al.} \cite{Franc}.

\subsection{Orthogonal Collocation}

Orthogonal collocation \footnotemark is a simultaneous (collocation) method that uses orthogonal polynomials to approximate the state and control functions. Orthogonal polynomials have several useful properties. The key concept is that a polynomial can be represented over some finite domain by its value at a special set of grid points over that domain. When represented in this form, it is easy to do fast and accurate numerical interpolation, differentiation, and integration of the polynomial \cite{Berrut2004}. 

\footnotetext{Orthogonal Collocation is often referred to as the {\em Pseudospectral Method} for transcription of an trajectory optimization problem. }

\subsection{H vs P scheme}

\par The simplest version of orthogonal collocation works by representing the entire trajectory as a single high-order orthogonal polynomial. The constraint on the dynamics (either integral or derivative form) is applied at the collocation points, which are chosen to be at the roots of the orthogonal polynomial. The cost function is then evaluated by numerical quadrature over the trajectory, which is just a weighted combination of the function's value at each collocation point. Convergence in this method is obtained by increasing the order of the polynomial (p-method). 

\par
\begin{figure}
\centering
\includegraphics[width=\textwidth]{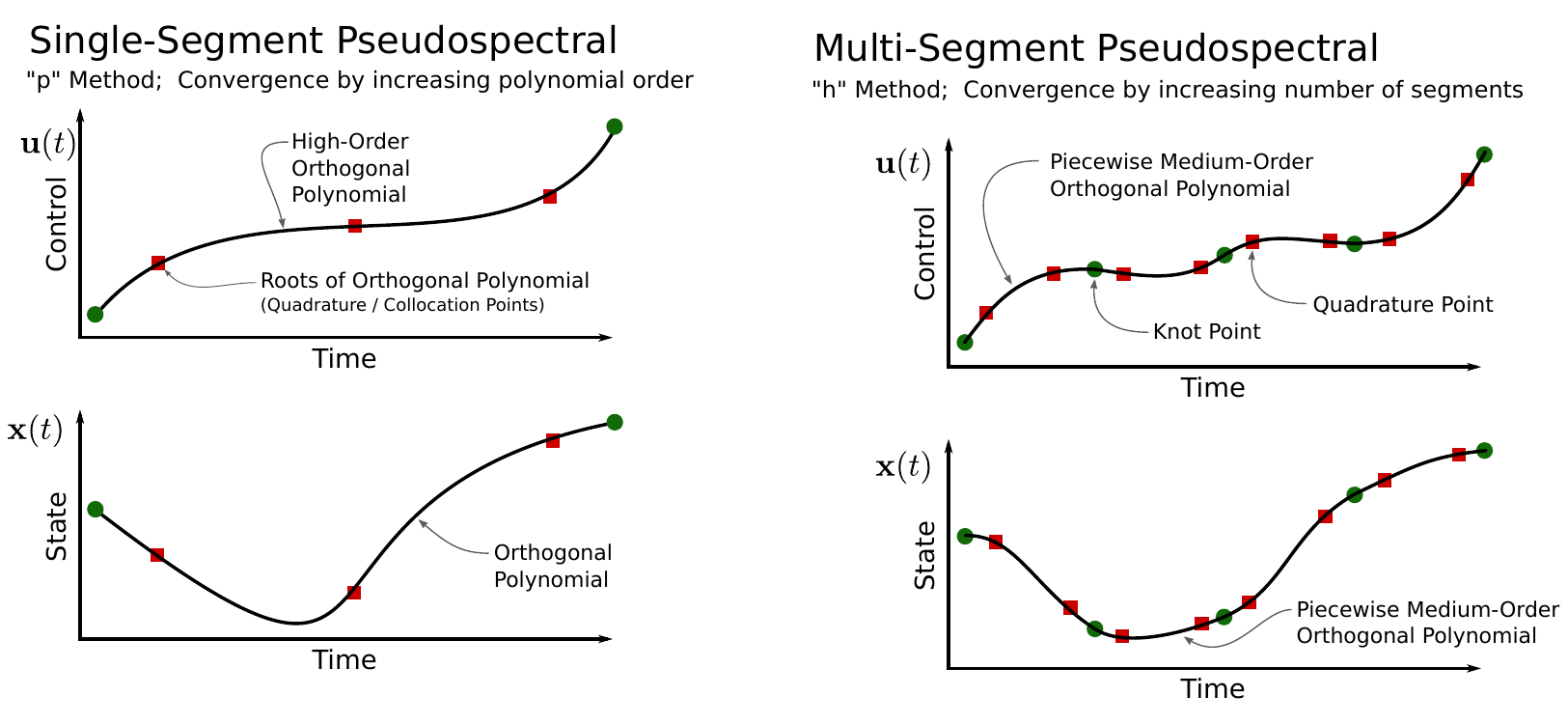}
\caption{Orthogonal Collocation Methods}
\label{fig:Transcription__Pseudospectral}
\end{figure}

 Representing the entire trajectory as a single orthogonal polynomial works well if the underlying solution is analytic. In many cases this is not true, such as when an actuator is saturated. In these cases the polynomial approximation will necessarily fail near the discontinuity. One solution is to represent the trajectory as a series of medium-order orthogonal polynomials. Each section of the trajectory is stitched together with defect constraints, like in multiple shooting, and the dynamics within a given segment are expressed using a constraint at the collocation points. The points where two segments are jointed are called {\em knot points}. Convergence in this type of method is obtained by increasing the number of trajectory segments (h-method). Figure \ref{fig:Transcription__Pseudospectral} shows a comparison of the single- and multi-segment orthogonal collocation methods.

\par The most sophisticated methods take this a step further, and adaptively refine the number and location of knot points along the trajectory and the order of the polynomial on each segment\cite{PattersonHagerRao}. This is the approach taken by the commercially available transcription algorithm GPOPS\cite{Patterson2013}.

\subsection{Special Cases}

Two commonly used simultaneous methods are {\em direct transcription} and {\em direct collocation}. Direct transcription is a simultaneous method that uses the integral form of the dynamics constraint. The control is represented as a piecewise-constant trajectory, and the state is piecewise linear. The decision variables in the optimization are the values for the control and state at each knot point.

\par
\begin{figure}
\centering
\includegraphics[width=\textwidth]{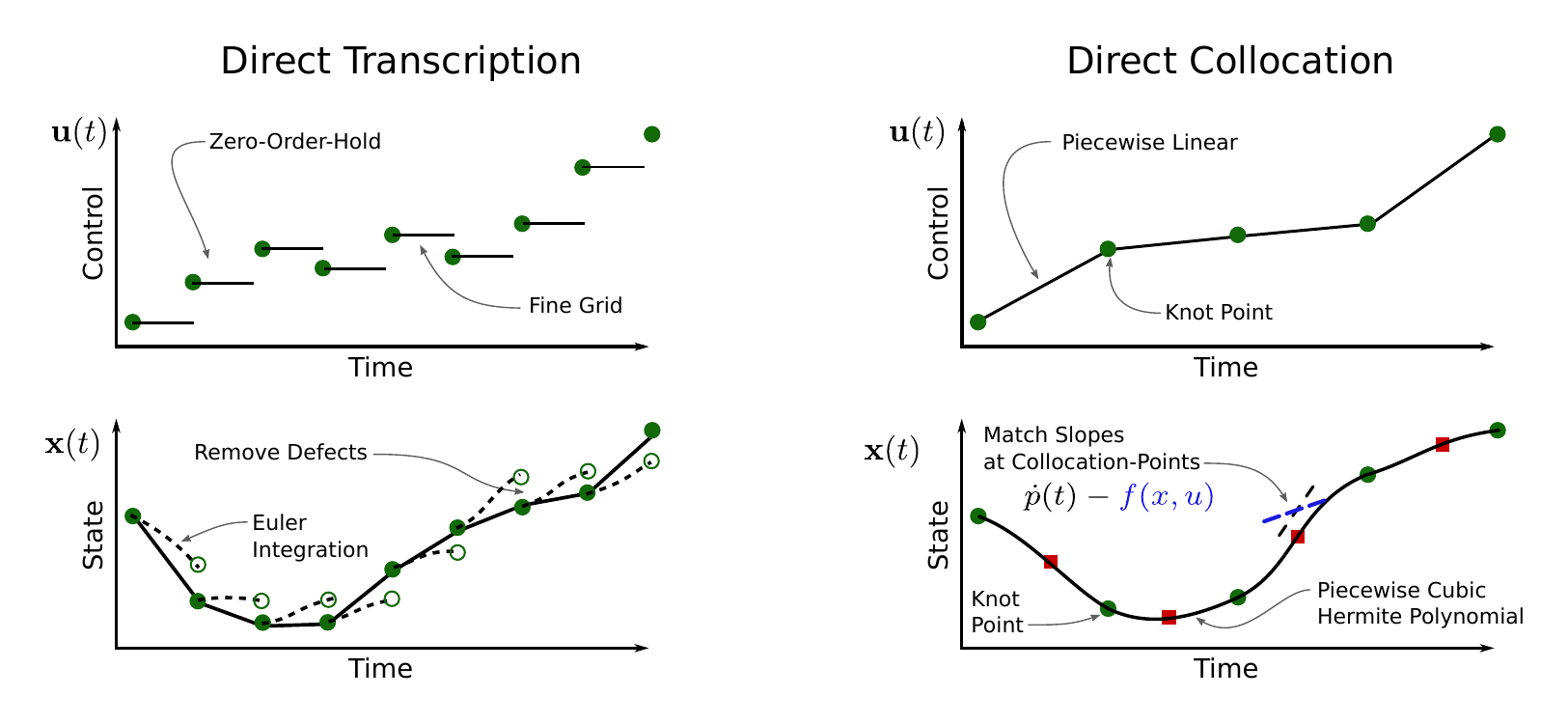}
\caption{Special Cases of simultaneous methods}
\label{fig:Transcription_Collocation}
\end{figure}

 Direct collocation is similar to direct transcription, but the input is represented as a piecewise-linear function of time, and the state is piecewise-cubic. The values of the state and control at each knot point are the decision variables. The slope of state is prescribed by the dynamics at each knot point. The collocation points are the mid-points of each cubic segment. The slope of the cubic at the collocation point is constrained to match the system dynamics at that point. Figure \ref{fig:Transcription_Collocation}.

%======================================
\section{Dealing with Hybrid Systems} \label{sec:HybridSystems}

A hybrid system is a dynamical system which has multiple phases of continuous dynamics, separated by discrete transitions.The simplest example of a hybrid system is a bouncing ball. It has a single continuous phase of motion (free fall), and a single discrete transition (collision with the ground). 

\par A more complicated hybrid system would be a walking robot. A few examples of continuous modes might be:

\begin{itemize} \setlength{\itemsep}{-4pt}%
    \item {\bf flight - } both feet in the air
    \item {\bf double stance - } both feet on the ground
    \item {\bf single stance - } one foot on the ground
\end{itemize}

There would then be a discrete transition between each of these modes. Some transitions do not alter the continuous state (such as lifting a foot off of the ground), while other transitions do change the continuous state (a foot colliding with the ground).

\par A naive way to handle such hybrid systems would be to bundle them up inside a big simulation or dynamics function and try to just make a simple transcription method work. This will almost certainly fail. The reason is that the underlying optimization algorithm is typically a smooth, gradient-based method (eg. SNOPT, IPOPT, FMINCON). The discrete transitions of the hybrid system cause the dynamics to be non-smooth. For example, a small change in one variable could cause the collision to happen at a different time, which would then cause a huge change in the defect constraints.

\par There are two generally accepted methods for dealing with hybrid systems: multi-phase methods, and through-contact methods. Multi-phase methods are faster and more accurate, but they require explicit knowledge of the sequence of transitions. Through-contact methods can deal with arbitrary sequences of contacts, but are slower and less accurate.

\par These two methods are analogous to the two different methods of simulation hybrid systems. Multi-phase transcription is like simulating a finite state machine using event detection. Through-contact is analogous to time-stepping.

\subsection{Multi-Phase Methods}

Multi-phase methods are a fairly simple extension of the basic multiple shooting or collocation algorithms. First, the user decides which sequence of phases should occur. For example, in a walking robot, this might be single stance (one foot on the ground) followed by double stance (two feet on the ground). Then the constraints are set up within each of those continuous phases exactly the same as in a simpler problem. Then a set of constraints are added to switch the continuous phases together, satisfying the transition equations.

\par One interesting side effect of this method is that the state of the system can be completely different in each phase of motion, provided that there is some sensible way to link them. This is useful, because it means that each set of continuous dynamics can be written in minimal coordinates.

\par The commercially available transcription program GPOPS II \cite{Patterson2013} provides a sophisticated interface for setting up and solving hybrid system trajectories using multi-phase methods.

\subsection{Through-Contact Methods}

Through-contact methods take a completely different approach. Instead of pushing the discontinuities in the dynamics to a special constraints between each phase, they directly handle the contact constraints at every grid point.

\par The key idea in through-contact optimization is that discontinuities are fine, provided that they are directly handled by constraints. This is done by writing the system dynamics in an impulse-based form, and leaving the contact impulses arbitrary (rather than algebraically eliminating them by assuming the contact mode). The (unknown) contact impulses are then treated as a control variable, which is subject to the following constraints at every grid point:

\begin{align*}
	&d_n > 0            \quad &\text{Contact seperation} \\
	&J_t \leq \mu J_n  \quad &\text{Contact force in friction cone} \\
	&d_n J_t = 0 \quad &\text{Contact force when touching}
\end{align*}

These constraints \footnote{These constraints form what is known as a {\em linear complementarity problem} (LCP).} will provide a unique solution for the contact impulses at every grid-point. This allows for an arbitrary sequence of contacts throughout the trajectory, which is particularly useful for complex behaviors, such as a biped crawling, or standing up from a laying down position, as shown by Mordatch et al \cite{Mordatch2012}.

\par Through-contact trajectory optimization will soon be included as part of Drake \cite{Tedrake2014}, a toolbox for doing robot control design and analysis. 

\subsection{Accuracy of each method}
\par If the sequence of contact modes is unknown, then a through contact method might be more accurate because it can find solution that would be excluded by the prescribed sequence of phases in the multi-phase method. That being said, if both methods manage to find roughly the same solution, then the multi-phase method will be far more accurate.

\par The multi-phase method can be incredibly accurate because the optimization is implicitly doing root finding for the transitions between each continuous phase of motion. In addition to that, the continuous trajectory can be easily represented by a high-order polynomial (either directly in collocation, or indirectly through a Runge-Kutta scheme in multiple shooting). These high-order trajectories are able to match the true solution very well, even with relatively large time steps, since the errors are proportional to the order of the polynomial.

\par The through-contact methods are inherently limited in their accuracy because the transitions between continuous phases of motion are constrained to grid-points, which are not able to be independently adjusted. This means that the approximation errors are proportional to the step-size, drastically limiting the ability to get high-accuracy solutions.

%=================================
\section{Setting up your problem}

\subsection{Grid Selection}

Both shooting and simultaneous methods have some free parameter(s) that control the grid that it is used to transcribe the problem. In general, you need a fine grid for an accurate solution, but this can cause convergence problems if the initial guess is not very good. A common solution is to use a very coarse grid to find an approximate solution, and then use this approximate solution as an initial guess for a second, or third, optimization using a finer grid. The commercially available program GPOPS II handles this grid-refinement automatically using a sophisticated algorithm \cite{Darby2011a}.

\subsection{Initialization}

Even a well-posed trajectory optimization is likely to fail with a poor initialization. One good method for initializing a trajectory is to guess a few points on the trajectory, and then fit a polynomial to these points. Then this polynomial can be differentiated once to get the first derivatives of the state, and again to get the joint accelerations. Inverse dynamics can then be used to compute the joint torques necessary to produce those accelerations. Then interpolate this rough guess to get the correct grid-points for either your multiple shooting or collocation method.

\par One thing that can go wrong with initialization is that you start the optimization with an infeasible, but locally optimal solution. This commonly occurs when you set the initial control functional to zero. This can be corrected by added a small random noise along the initial guess at the actuator. This has the added benefit of starting you optimization from a new place on each iteration, which can sometimes be helpful to detect local minima.

\par Sometimes the optimization will fail, even with a reasonable initialization, if the cost function is too complicated\footnote{Cost of transport, which is total energy used by a robot divided by it's weight multiplied by the distance travelled, is a particularly difficult cost function to optimize.}. One solution is to write the optimization to use an alternate cost function, such as torque-squared, that is very simple to optimize. Use this alternate cost function when running the first, coarse grid, optimization, and then once a feasible trajectory is found, change to the more difficult cost function. Sometimes this can cause problems with local minima, but it is a good first thing to try when debugging.

\subsection{Local Minimum}

\par
\begin{figure}
\centering
\includegraphics[width = 0.5\textwidth]{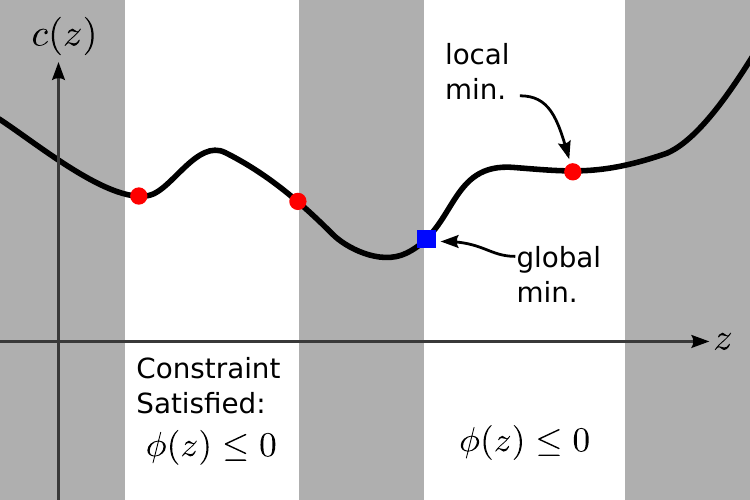}
\caption{Cartoon of a local minimum in a simple 1D constrained optimization problem. Figure inspired by Russ Tedrake's lecture notes for his Underactuated Robotics class.}
\label{fig:localMinCartoon}
\end{figure}

 Once you have a `converged' solution to the trajectory optimization problem, it is good practice to check if the solution is the 'real' (global) solution. There is no practical way to prove that you have found the global solution in most cases, but you can run the optimization from several initial guesses and check that they all converge to the same solution. Figure \ref{fig:localMinCartoon} shows a cartoon of a one-dimensional constrained optimization problem with three sub-optimal local minima and a single global minimum.

%==========================
\section{Smoothness and Consistency}

In order for a trajectory optimization program to converge, all of the internal function calls must be smooth and consistent:

\begin{itemize}
	\item {\bf Smooth - } The output of the function is continuous, and at least twice differentiable at all points. Ideally the derivatives should all remain small - sharp (continuous) corners will still cause problems.
	\item {\bf Consistent - } The output of the function is repeatable and smooth across multiple calls; the computer executes the same exact lines of code on every call (no \texttt{if} statements).
\end{itemize}

These two requirements sound simple, but in practice they are remarkably difficult to deal with. In this section I will attempt to cover several common places that errors are made. Note that these requirements apply to both the constraints and objective function.

\subsection{Dynamics}

Walking robots are hybrid systems - they have both continuous and discrete aspects to their dynamics. These discrete dynamics must be carefully handled, either by using through-contact or multi-phase transcription methods, as discussed in Section \ref{sec:HybridSystems}.

Another common source of error is in the simulation of the dynamics function for multiple shooting methods. In shooting methods, it is important to use an explicit, fixed-step integration algorithm, rather than an adaptive method like \texttt{ode45}.

\par Sometimes it might be tempting to put some sort of noise into the dynamics. This guarantees that the function will be inconsistent, and will almost certainly cause the optimization program to fail. In general, there is no good reason to put a call to a random number generator inside of a trajectory optimization program.

\subsection{Objective Function}

There are many sources of discontinuity inside of an objective function, typically due to seeming innocuous functions. A few examples are taking the maximum or minimum of a list of numbers, the absolute value function, clamp function\footnote{Also called saturation}, and the ramp function.

\par The best way to deal with these discontinuous functions is by using a constraint to enforce the discontinuity. This works because the non-linear program solver has special tools build-in for dealing with constraints. Equation \ref{eq:AbsVal} shows the correct way to handle the function $|x|$. Other examples are provided in \cite{JohnT.Betts2001}.

\begin{equation} \label{eq:AbsVal}
|x| = x_1 + x_2 \quad \text{subject to} \quad

  \begin{cases}
	x=x_1-x_2 \\
	x_1>=0 \\ 
	x_2>=0	
  \end{cases}
\end{equation}

\par Another approach to dealing with discontinuous functions is to use `smooth' approximations, as shown below:

\begin{align*}
|x| \quad \approx \quad & x \tanh \left(x / \alpha \right) \\
\max(\mathbf{x}) \quad \approx \quad &  \alpha \log\left(\sum \exp \frac{\mathbf{x}}{\alpha} \right)
\end{align*}

\par This method works, but it will typically be less accurate and take longer to converge than the constraint method. Even though the smoothed functions are continuous, there is a tradeoff between accuracy and convergence. If the smoothing is large, then the optimization program will converge, but with an inaccurate solution. If the smoothing is small, the program will converge slowly, or perhaps not at all. One solution is to iteratively reduce the smoothing, starting with heavy smoothing to get a feasible solution, and then reducing the smoothing to zero in on the answer. 

\par There are two approaches to making smooth approximations of functions. The first is to use exponential functions that asymptotically approach the true function far from the discontinuity. The second is to use a piecewise version of the function, where the region near the discontinuity is replaced with a polynomial approximation.

%======================================
\section{Software Specifics}

Once a trajectory optimization problem been transcribed by either multiple shooting or collocation, it must then be solved by nonlinear constrained optimization solver. Two of the most popular algorithms are SNOPT and FMINCON. The following section describe some specific details that are helpful when working with each of these programs.

\par It seems that FMINCON solves the multiple shooting problem by first finding a feasible solution, and then attempting to optimize it. This means that FMINCON can handle a worse initial guess than SNOPT, but it is a little worse at finding the true optimal solution, because it does not allow for much flexibility in the constraints. 

\par Instead of writing your own transcription algorithm, it is possible to buy a commercial version. The best available software now is probably GPOPS II, which internally calls either SNOPT or IPOPT (another solver, not discussed here). 

%%~~~~~~~~~~~~~~~~~~~~~~~~~~~~~~~~~~~~~~~~~~~~~~~~~~~~~~~~~~~~~~~~~~
\subsection{SNOPT \cite{Gill2005}}

	\par A single function call is used for the objective function and all constraints. This function returns a vector, and the first element of the vector is the value of the objective function. All of the following elements are constraints. SNOPT does not require the user to specify which constraints are linear or nonlinear - it automatically computes this.

	\par Since SNOPT combines the objective function and the constraint function, you always need to provide a value for the objective function, even if you are only solving a feasibility problem. In this case, the value of the objective function should be set to zero. More importantly, you must also specify that the bounds of the objective function are also zero. In code this looks like: \\
					\texttt{F(1)==0, Flow(1)==0, Fupp(1)==0}

	\par Do not use a constant term anywhere in your objective function or constraints. All constant terms must be moved to the boundary vectors (\texttt{Flow, Fupp}). This is because of how SNOPT internally stores and estimates the gradients - it will numerically remove any constant terms.

	\par SNOPT does not allow the user to pass any parameters to the objective function, at least when called from Matlab. One way to get around this is to create a single struct of parameters, and save it as a global variable. Once this is done, you can then just read the parameters off the global variable.

%%~~~~~~~~~~~~~~~~~~~~~~~~~~~~~~~~~~~~~~~~~~~~~~~~~~~~~~~~~~~~~~~~~~
\subsection{FMINCON \cite{MatlabOptimizationToolbox2014}}

FMINCON requires that the objective function and constraint function are in separate functions. This creates a problem, because both the constraints and objective need to evaluate the dynamics at every grid-point. To avoid doing all of the work twice, it is best to use a persistent variable (this would be a static variable in C++). Create a single function that does the integration of the dynamics, and have it store the last input state and output solution. When it is called, it checks if the state is the same as last time - if it is, then just return the previous solution.

\par There is a bug of sorts in FMINCON, which prevents you from using state bounds to constrain a state to a specific value. Basically, saying: $0<x<=0$ causes FMINCON to try to divide by a step size of machine precision when doing its finite differencing. Much better to place these sorts of things as an equality constraint, which are handled in a special way.

%%~~~~~~~~~~~~~~~~~~~~~~~~~~~~~~~~~~~~~~~~~~~~~~~~~~~~~~~~~~~~~~~~~~
\subsection{GPOPS2 \cite{Patterson2013}}

GPOPS2 is a commercially available software which implements {\em ``a general-pupose software for solving optimal control problems using variable-order adaptive orthogonal collocation methods together with sparse nonlinear programming''}

This program is excellent at solving trajectory optimization problems which have a known sequence of contact modes, making an extremely general solver. In my experience it is generally much faster and more accurate than using a more simple multiple shooting or collocation method.

\section{Miscellaneous}

This section covers a few odd topics that are useful to understand, but do not fit elsewhere.

\subsection{Regularization}
It is possible to create an optimal trajectory problem which does not have a unique solution. This will generally cause convergence problems, as the optimization program bounces between seemingly equivalent solutions. This problem is solved by adding a small regularization term to the objective function, which forces a unique solution. For dynamical systems, I have found that adding a small input-squared term to the cost function is generally sufficient. I've found regularization terms that are 6-8 orders of magnitude smaller than the primary objective term are often still effective.

\subsection{Constraint on Controls}
The correct way to apply a non-linear constraint to a control at the boundary of a trajectory is to create a dummy state to represent the control, and then let the optimization program determine the {\em derivative} of the control. The system dynamics can then be used to ensure feasibility of the control and its derivative. This technique is particularly useful for enforcing that the contact forces of a walking robot stay within their friction cone.

\subsection{Optimizing a Parameter}
Suppose that you would like to find an optimal trajectory, but there is at least one free parameter for your design. It is tempting to just add this parameter as an additional term when passed to the underlying optimization function. This is a bad idea, because it couples the Jacobian of the constraints, which is almost as bad as solving the problem using single shooting. The correct thing to do is to add a control to the system, and use the control instead of the parameter. A special constraint is then added to ensure that the control remains constant throughout the trajectory. This will quickly solve for the best possible parameter choice, while keeping the Jacobian of the constraints and objective function sparse.

%=================================================

\section{Hammer Example}

I've created a simple example to demonstrate a simple multiple shooting transcription algorithm (using Matlab's FMINCON to solve the underlying optimization). In this example, the goal is to find a trajectory for hammer that is periodically striking a surface. The hammer is modeled as a point-mass pendulum and powered by a torque source. The cost of using the torque source is modeled as the integral of torque squared against time.

\par A series of plots are shown to detail the progress that the optimization algorithm makes towards the solution. On the left side of each figure is a plot of the state-space trajectory. Note that there is a discrete jump in the trajectory when the hammer strikes the surface. Below this trajectory is a plot that shows how the cost function accumulates over time. On the right side of each figure is a set of three plots, each one showing a single component of the trajectory against time (hammer angle, hammer angular rate, and torque applied to hammer). The last plot has a marker for every gridpoint in the state-space trajectory, and shows the final number of iterations. 

\par In the first plot the trajectories are jagged and have red and black lines. The red lines show the {\it defects} in the trajectory - as the optimization runs these become arbitrarily small. 

\par The full code for running this example is available on my website: \\
\footnotesize{\texttt{http://ruina.tam.cornell.edu/research/MatthewKelly}}

\FloatBarrier

\onecolumn

\begin{figure} 
    \centering 
    \includegraphics[width = \columnwidth]{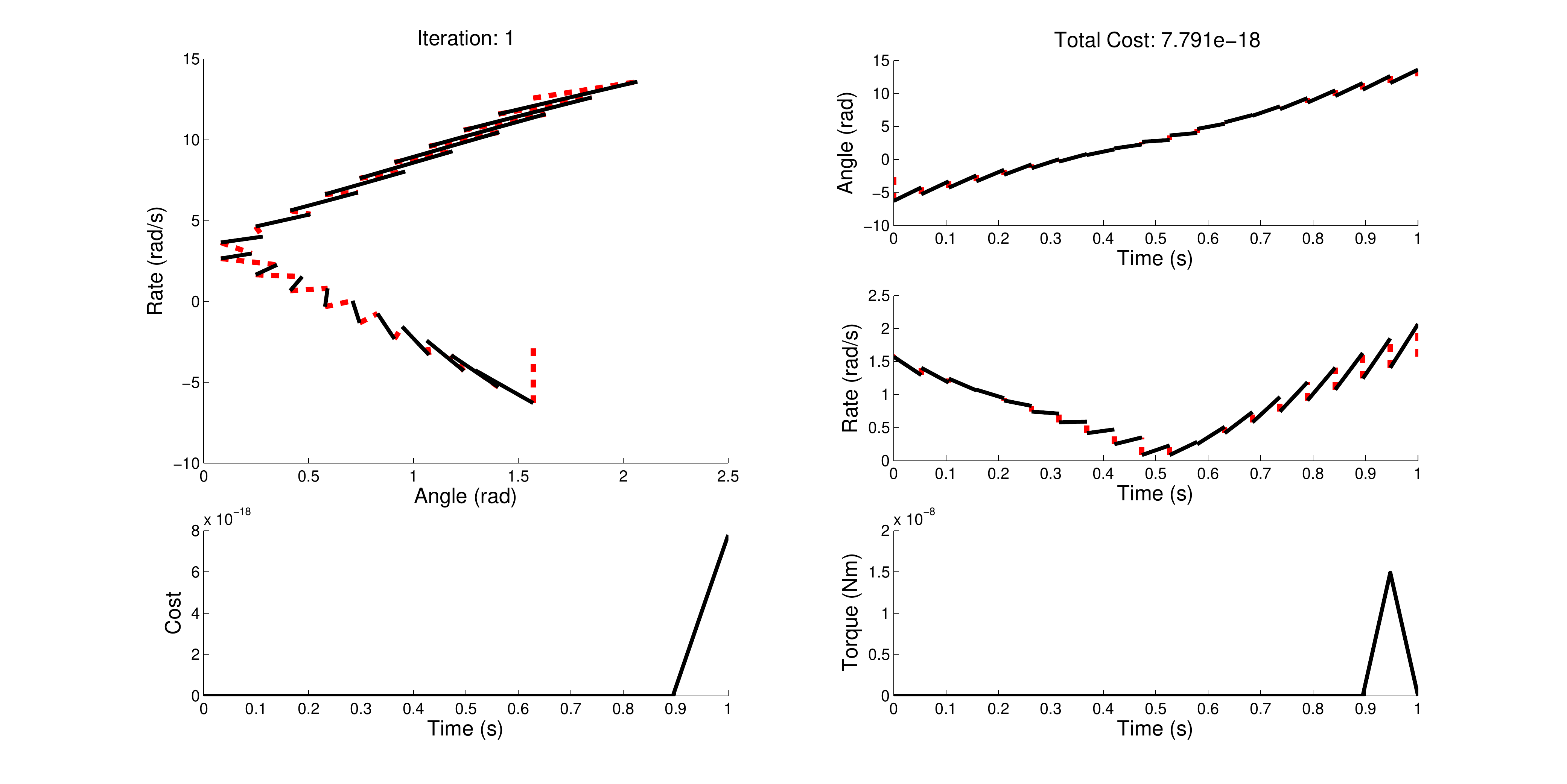}    
    \caption{The first iteration of the optimization program. Notice that FMINCON has added a small perturbation to the guess at the control, to help numerically calculate the Jacobian. The initial guess was a trajectory made of two straight lines, and the initial defects are clearly visible.}     
    \label{fig: MS_Fig_Iter_1}     
\end{figure} 

%\begin{figure} 
%    \centering 
 %   \includegraphics[width = \columnwidth]{img/MS_Fig_Iter_2.pdf}    
 %   \caption{}     
  %  \label{fig: MS_Fig_Iter_2}     
%\end{figure} 

\begin{figure} 
    \centering 
    \includegraphics[width = \columnwidth]{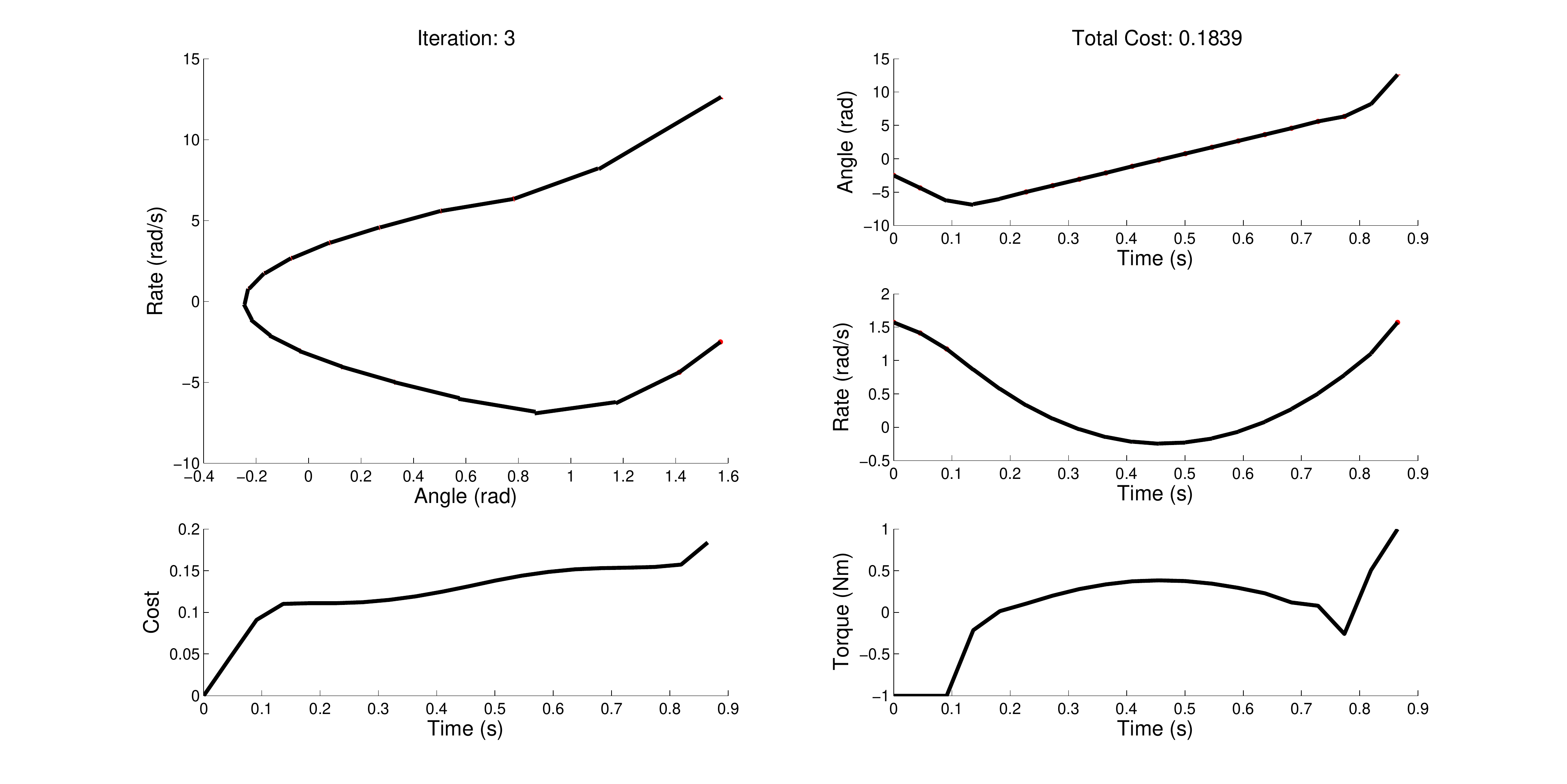}    
    \caption{Even after a few iterations, the defects are no longer visible on the plot - the trajectory is nearly feasible . Notice that the optimization has a first guess at the actuation, and a non-zero cost function.}     
    \label{fig: MS_Fig_Iter_3}     
\end{figure} 

%\begin{figure} 
%    \centering %
%    \includegraphics[width = \columnwidth]{img/MS_Fig_Iter_5.pdf}    
   % \caption{}     
 %   \label{fig: MS_Fig_Iter_5}     
%\end{figure} 

\begin{figure} 
    \centering 
    \includegraphics[width = \columnwidth]{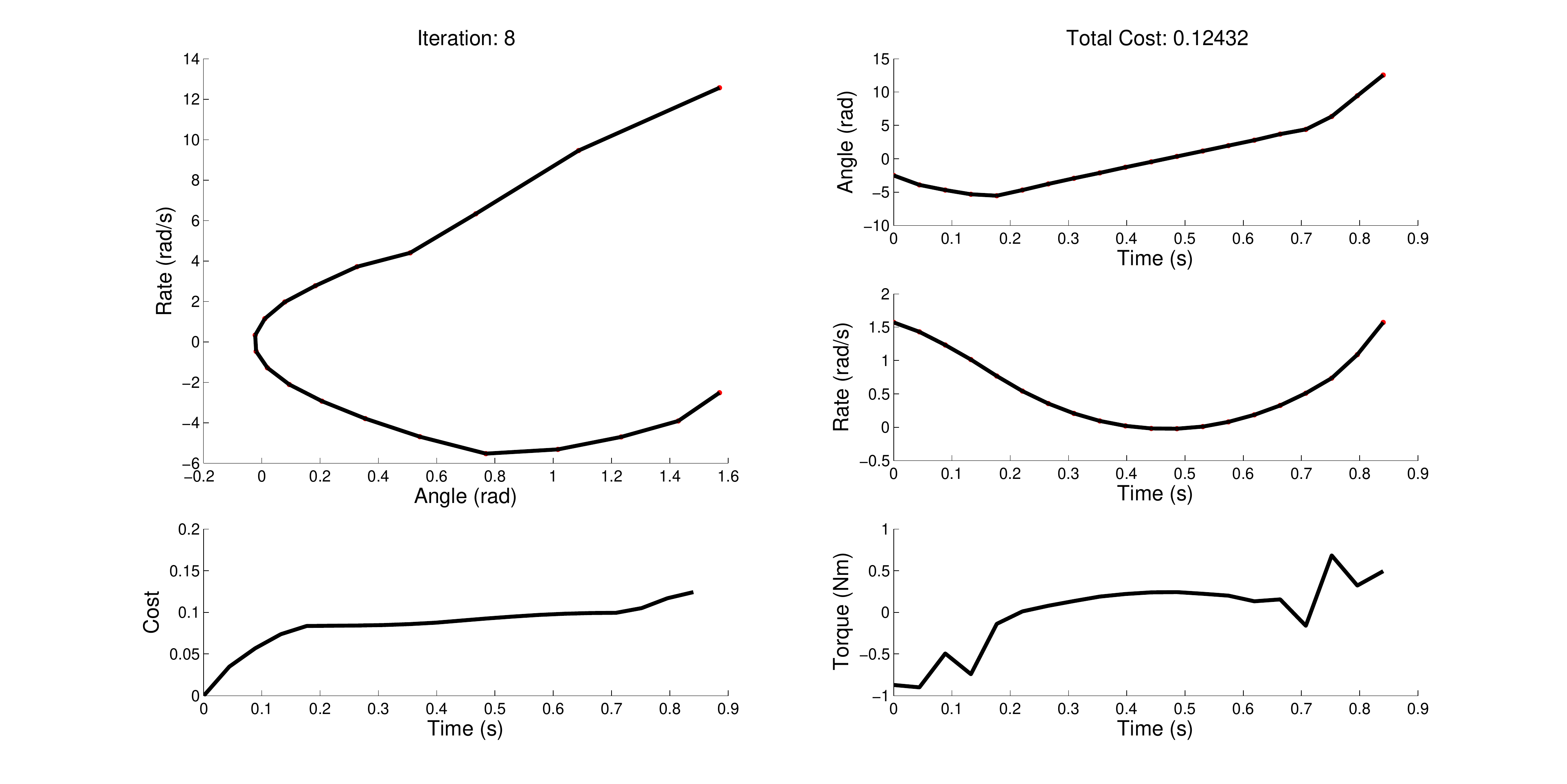}    
    \caption{Now the trajectory has roughly stabilized, and the optimization program is trying to make small changes to the control to reduce the total cost of the trajectory. Interestingly, the trajectory is smooth, but the control is still fairly discontinuous. This is indicative of an solution that has not fully converged.}     
    \label{fig: MS_Fig_Iter_8}     
\end{figure} 

%\begin{figure} 
%    \centering 
%    \includegraphics[width = \columnwidth]{img/MS_Fig_Iter_13.pdf}    
 %   \caption{}     
 %   \label{fig: MS_Fig_Iter_13}     
%\end{figure} 

%\begin{figure} 
%    \centering 
%    \includegraphics[width = \columnwidth]{img/MS_Fig_Iter_21.pdf}    
%    \caption{Pretty similar to the previous iteration, but with a slightly smaller cost and smoother control function.}     
%    \label{fig: MS_Fig_Iter_21}     
%\end{figure} 

%\begin{figure} 
 %   \centering 
 %   \includegraphics[width = \columnwidth]{img/MS_Fig_Iter_34.pdf}    
 %   \caption{}     
  %  \label{fig: MS_Fig_Iter_34}     
%\end{figure} 

%\begin{figure} 
%    \centering 
%    \includegraphics[width = \columnwidth]{img/MS_Fig_Iter_55.pdf}    
%    \caption{ALmost converged, but there are still a few wiggles to get out of the control solution.}     
%    \label{fig: MS_Fig_Iter_55}     
%\end{figure} 

%\begin{figure} 
%    \centering 
%    \includegraphics[width = \columnwidth]{img/MS_Fig_Iter_89.pdf}    
 %   \caption{}     
 %   \label{fig: MS_Fig_Iter_89}     
%\end{figure} 

%\begin{figure} 
%    \centering 
 %   \includegraphics[width = \columnwidth]{img/MS_Fig_Iter_144.pdf}    
 %   \caption{}     
 %   \label{fig: MS_Fig_Iter_144}     
%\end{figure} 

\begin{figure} 
    \centering 
    \includegraphics[width = \columnwidth]{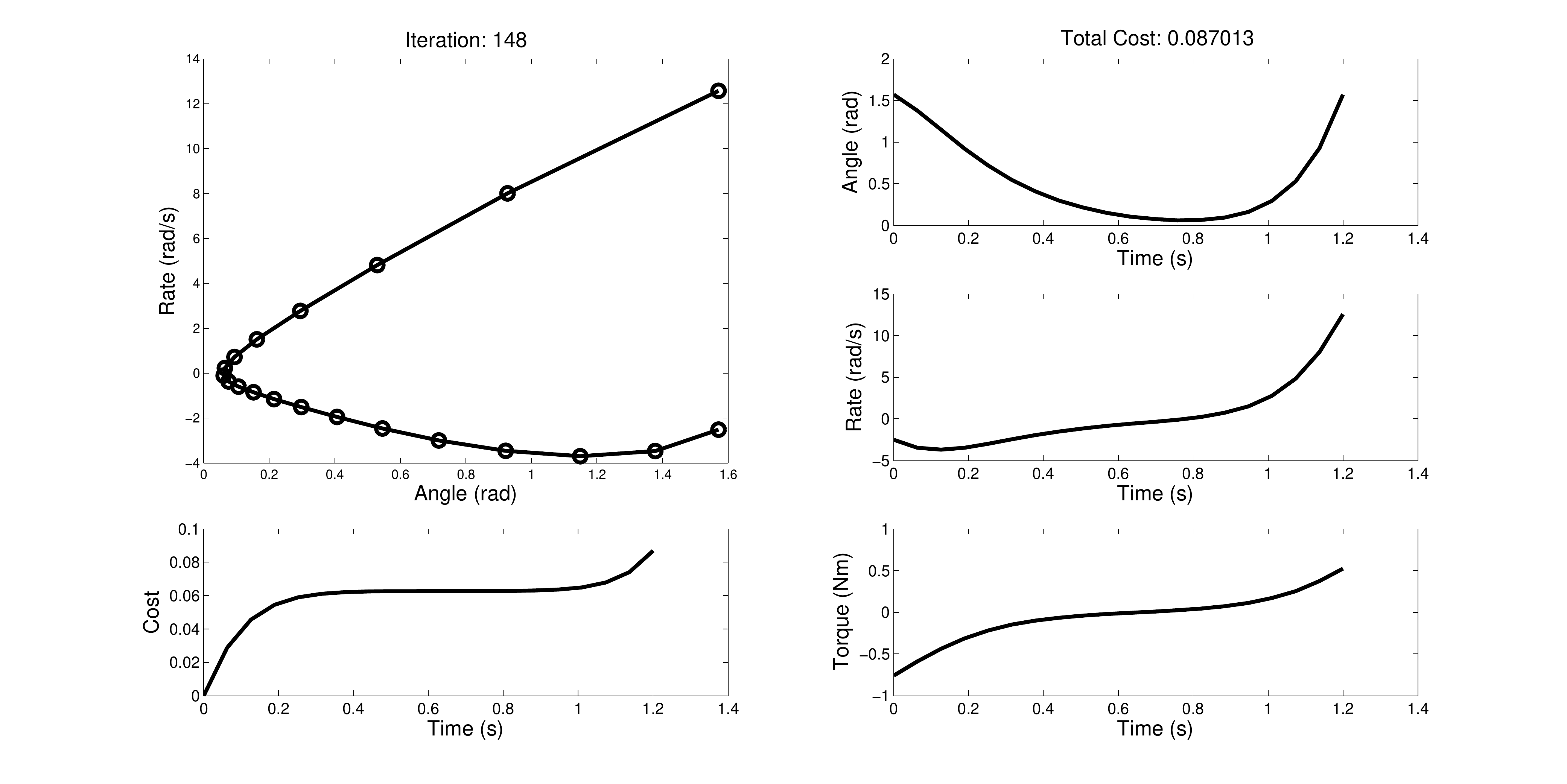}    
    \caption{ After 148 iterations the optimization program terminated. Notice that all components of the trajectory are smooth, including the control function. I've put small circles over the start and end of each segment of the trajectory, which now line up almost perfectly.}     
    \label{fig: MS_Fig_Iter_148}     
\end{figure}

\FloatBarrier

\bibliographystyle{plain}
\bibliography{../../library}

\end{document}